\documentclass[10pt]{article}
\usepackage{latexsym}
\usepackage{amsfonts}
\usepackage{enumerate}
\usepackage{multicol}
\usepackage{graphicx}
\usepackage{amssymb}
\usepackage{amsmath}
\usepackage{epic}

\title{The Minimum Size of Signed Sumsets   \\[.4in]}

\author{B\'{e}la Bajnok\footnote{Corresponding author} \\[.1in] {\small Department of Mathematics, Gettysburg College} \\
{\small 300 N. Washington Street, Gettysburg, PA 17325-1486 USA} \\{\small E-mail:  bbajnok@gettysburg.edu} \\ [.2in]
and \\[.2in]
Ryan Matzke \\[.1in] {\small Department of Mathematics, Gettysburg College} \\
{\small 300 N. Washington Street, Gettysburg, PA 17325-1486 USA} \\{\small E-mail:  matzry01@gettysburg.edu} \\ [.2in]
 \\[.4in]}

\date{December 1, 2014}

\newtheorem{thm}{Theorem}

\newtheorem{cor}[thm]{Corollary}
\newtheorem{prop}[thm]{Proposition}
\newtheorem{conj}[thm]{Conjecture}

\begin{document}

\maketitle

\begin{abstract}

For a finite abelian group $G$ and positive integers $m$ and $h$, we let 
$$\rho(G, m, h) = \min \{ |hA| \; : \; A \subseteq G, |A|=m\}$$ and
$$\rho_{\pm} (G, m, h) = \min \{ |h_{\pm} A| \; : \; A \subseteq G, |A|=m\},$$ where $hA$ and $h_{\pm} A$ denote the $h$-fold sumset and the $h$-fold signed sumset of $A$, respectively.  The study of $\rho(G, m, h)$ has a 200-year-old history and is now known for all $G$, $m$, and $h$.  Here we prove that $\rho_{\pm}(G, m, h)$ equals $\rho (G, m, h)$ when $G$ is cyclic, and establish an upper bound for $\rho_{\pm} (G, m, h)$ that we believe gives the exact value for all $G$, $m$, and $h$.  
      
\end{abstract}



\thispagestyle{empty}

\section{Introduction}

Let $G$ be a finite abelian group written with additive notation.  For a nonnegative integer $h$ and a nonempty subset $A$ of $G$, we let $hA$ and $h_{\pm}A$ denote the $h$-fold {\em sumset} and the $h$-fold {\em signed sumset} of $A$, respectively; that is, for an $m$-subset $A=\{a_1, \dots, a_m\}$ of $G$, we let
$$hA=\{ \Sigma_{i=1}^m \lambda_i a_i \; : \; (\lambda_1,\dots,\lambda_m) \in \mathbb{N}_0^m, \; \Sigma_{i=1}^m \lambda_i=h\}$$ and 
$$h_{\pm} A=\{ \Sigma_{i=1}^m \lambda_i a_i \; : \; (\lambda_1,\dots,\lambda_m) \in \mathbb{Z}^m, \; \Sigma_{i=1}^m |\lambda_i|=h\}.$$  

While signed sumsets are less well-studied in the literature than sumsets are, they come up naturally: For example, in \cite{BajRuz:2003a}, the first author and Ruzsa investigated the {\em independence number} of a subset $A$ of $G$, defined as the maximum value of $t \in \mathbb{N}$ for which $$0 \not \in \cup_{h=1}^t h_{\pm}A$$ (see also \cite{Baj:2000a} and \cite{Baj:2004a}); and in \cite{KloLev:2003a}, Klopsch and Lev discussed the {\em diameter} of $G$ with respect to $A$, defined as the minimum value of $s \in \mathbb{N}$ for which $$\cup_{h=0}^s h_{\pm}A=G$$ (see also \cite{KloLev:2009a}).  The independence number of $A$ in $G$ quantifies the ``degree'' to which $A$ is linearly independent in $G$ (no subset is ``completely'' independent), while the diameter of $G$ with respect to $A$ measures how ``effectively'' $A$ generates $G$ (if at all).  Note that $h_{\pm}A$ is always contained in $h(A \cup -A)$, but this may be a proper containment when $h \geq 2$.      

For a positive integer $m \leq |G|$, we 
let $$\rho(G, m, h) = \min \{ |hA| \; : \; A \subseteq G, |A|=m\}$$ and   
$$\rho_{\pm}(G, m, h) = \min \{ |h_{\pm}A|  \; : \; A \subseteq G, |A|=m\}$$ (as usual, $|S|$ denotes the size of the finite set $S$).
The value of $\rho(G, m, h)$ has a long and distinguished history and has been determined for all $G$, $m$, and $h$; in this paper we attempt to find $\rho_{\pm}(G, m, h)$.   

We start by a brief review of the case of sumsets.  In 1813, for prime values of $p$, Cauchy  \cite{Cau:1813a} found the minimum possible size of
$$A+B=\{a+b \; : \; a \in A,\; b \in B \}$$ among subsets $A$ and $B$ of given sizes in the cyclic group $\mathbb{Z}_p$.  In 1935, Davenport \cite{Dav:1935a} rediscovered Cauchy's result, which is now known as the Cauchy--Davenport Theorem.  (Davenport was unaware of Cauchy's work until twelve years later; see \cite{Dav:1947a}.)

\begin{thm}[Cauchy--Davenport Theorem] \label{Cauchy--Davenport}
If $A$ and $B$ are nonempty subsets of the group $\mathbb{Z}_p$ of prime order $p$, then 
$$|A+B| \geq \min \{p, |A|+|B|-1\}.$$
\end{thm}
It can easily be seen that the bound is tight for all values of $|A|$ and $|B|$, and thus
$$ \rho (\mathbb{Z}_p, m, 2)=\min\{p,2m-1\}.$$

After various partial results, the general case was finally solved in 2006 by Plagne \cite{Pla:2006a} (see also \cite{Pla:2003a}, \cite{EliKer:2007a}, and  \cite{EliKerPla:2003a}).  To state the result, we introduce the function
$$u(n,m,h)=\min \{f_d (m,h) \; : \; d \in D(n)\},$$ where $n$, $m$, and $h$ are positive integers, $D(n)$ is the set of positive divisors of $n$, and
$$f_d(m,h)=\left(h\left \lceil m/h \right \rceil-h +1 \right) \cdot d.$$  
(Here $u(n,m,h)$ is a relative of the Hopf--Stiefel function used also in topology and bilinear algebra; see, for example, \cite{EliKer:2005a}, \cite{Kar:2006a}, \cite{Pla:2003a}, and  \cite{Sha:1984a}.)

\begin{thm} [Plagne; cf.~\cite{Pla:2006a}] \label{value of u}
Let $n$, $m$, and $h$ be positive integers with $m \leq n$.  For any abelian group $G$ of order $n$ we have
$$\rho (G, m, h)=u(n,m,h).$$
\end{thm}

Turning now to $\rho_{\pm} (G, m, h)$, we start by observing that $$\rho_{\pm} (G,m,0)=1$$ and $$\rho_{\pm} (G,m,1)=m$$ for all $G$ and $m$.  To see the latter equality, it suffices to verify that one can always find a {\em symmetric} subset of size $m$ in $G$, that is, an $m$-subset $A$ of $G$ for which $A=-A$.   Therefore, from now on, we assume that $h \geq 2$.

We must admit that our study of $\rho_{\pm} (G, m, h)$ resulted in quite a few surprises.  For a start, we noticed that, in spite of the fact that $h_{\pm}A$ is usually much larger than $hA$ is, the equality $$\rho_{\pm} (G, m,h)=\rho (G, m,h)$$ holds quite often; it is an easy exercise to verify that, among groups of order 24 or less, equality holds with only one exception: $\rho_{\pm} (\mathbb{Z}_3^2, 4,2)=8$ while $\rho (\mathbb{Z}_3^2, 4,2)=7.$  In fact, we can prove that $\rho_{\pm} (G, m,h)$ agrees with $\rho (G, m,h)$ for all cyclic groups $G$ and all $m$ and $h$ (see Theorem \ref{cyclic} below).

However, in contrast to $\rho (G, m,h)$, the value of $\rho_{\pm} (G, m,h)$ depends on the structure of $G$ rather than just the order $n$ of $G$.  Suppose that $G$ is of type $(n_1,\dots,n_r)$, that is,   
$$G \cong \mathbb{Z}_{n_1} \times \cdots \times \mathbb{Z}_{n_r},$$ where $n_1 \geq 2$ and $n_i$ divides $n_{i+1}$ for each $i \in \{1,\dots, r-1\}$.  We exhibit a specific subset $D(G,m)$ of $D(n)$ with which the quantity
$$u_{\pm}(G,m,h)= \min \{f_d(m,h) \; : \; d \in D(G,m) \}$$ provides an upper bound for $\rho_{\pm} (G, m,h)$ (see Theorem \ref{u pm with f} below).  Therefore, to get lower and upper bounds for $\rho_{\pm} (G, m,h)$, we minimize $f_d(m,h)$ for all $d \in D(n)$ and for $d \in D(G,m)$, respectively:
$$\min\{ f_d(m,h) \; : \; d \in D(n)\} \leq \rho_{\pm} (G, m,h) \leq \min \{f_d(m,h) \; : \; d \in D(G,m) \}.$$
In fact, we also conjecture that $$\rho_{\pm} (G, m,h)=u_{\pm}(G,m,h)$$ holds in all but one very special situation (see Conjecture \ref{conj for rho pm} below).  

Further surprises come from the inverse problem of trying to classify subsets that yield the minimum signed sumset size.  To start with, we point out that it is not always symmetric sets that work best.  As an example, consider $\rho_{\pm} (\mathbb{Z}_5^2, 9, 2)$.  
One can see that for any 9 elements of $\pm a +H$, where $H$ is any subgroup of size 5 and $a \not \in H$, we have $$2_{\pm}A=H \cup (\pm 2a +H),$$ so $$\rho_{\pm} (\mathbb{Z}_5^2, 9, 2)=\rho(\mathbb{Z}_5^2, 9, 2)=15.$$  Here $A$ is not symmetric but is {\em near-symmetric}: it becomes symmetric once one of its elements is removed.  However, we can verify that for any symmetric subset $A$ of size 9,  $2_{\pm}A$ must have size 17 or more, as follows:   
If A contains a subgroup $H$ of size 5, then with any $a \in A \setminus H$, the 2-fold signed sumset of $A$ will contain the 17 distinct elements of $H$, $\pm a+H$, and $ \{\pm 2a\};$  while if $A$ contains no subgroup of size 5, then $$A \cap \{2a \; : \; a \in A\} = \{0\},$$ so $$|2_{\pm}A| \geq |A|+|\{2a \; : \; a \in A\}|-1=17.$$

And that's not all: sometimes it is best to take an {\em asymmetric} set, a set $A$ where $A$ and $-A$ are disjoint.  It is easy to check that, in the example of $\rho_{\pm}(\mathbb{Z}_3^2,4,2)=8$ mentioned above, with a 4-subset $A$ of $\mathbb{Z}_3^2$ we get $2_{\pm}A=\mathbb{Z}_3^2 \setminus \{0\}$ when $A$ is asymmetric, and $2_{\pm}A=\mathbb{Z}_3^2$ in all other cases.

We have thus seen that sets that minimize signed sumset size may be symmetric, near-symmetric, or asymmetric---we can prove, however, that there is always a set that is of one of these three types
 (see Theorem \ref{symmetry thm} below).

With this paper we aim to introduce the question of finding the minimum size of signed sumsets.  Our approach here is entirely elementary.  In the follow-up paper \cite{BajMat:2014b}, we investigate the question in elementary abelian groups, where, using deeper results from additive combinatorics, we are able to assert more.

\section{The role of symmetry}

Given a group $G$ and a positive integer $m \leq |G|$, we define a certain collection ${\cal A}(G,m)$ of $m$-subsets of $G$.  
We let 
\begin{itemize}
  \item $\mathrm{Sym}(G,m)$ be the collection of {\em symmetric} $m$-subsets of $G$, that is, $m$-subsets $A$ of $G$ for which $A=-A$;
  \item $\mathrm{Nsym}(G,m)$ be the collection of {\em near-symmetric} $m$-subsets of $G$, that is, $m$-subsets $A$ of $G$ that are not symmetric, but for which $A\setminus \{a\}$ is symmetric for some $a \in A$; 
  \item $\mathrm{Asym}(G,m)$ be the collection of {\em asymmetric} $m$-subsets of $G$, that is, $m$-subsets $A$ of $G$ for which $A \cap (-A)=\emptyset$.  
  \end{itemize}  We then let 
$${\cal A}(G,m)=\mathrm{Sym}(G,m) \cup \mathrm{Nsym}(G,m)\cup \mathrm{Asym}(G,m).$$  In other words, ${\cal A}(G,m)$ consists of those $m$-subsets of $G$ that have exactly $m$, $m-1$, or $0$ elements whose inverse is also in the set. 

\begin{thm} \label{symmetry thm}
For every $G$, $m$, and $h$, we have
$$\rho_{\pm} (G,m,h)= \min \{|h_{\pm} A| \; : \; A \in {\cal A}(G,m)\}.$$

\end{thm}

{\em Proof:}  Since our claim is trivial when $m \leq 2$, we assume that $m \geq 3$. 

For a subset $S$ of $G$, let us define its {\em degree of symmetry}, denoted by $\mathrm{sdeg}(S)$, as the number of elements of $S$ that are also elements of $-S$.   We shall prove that for any $m$-subset $B$ of $G$ with $$1 \leq \mathrm{sdeg}(B) \leq m-2,$$ there is an $m$-subset $B'$ of $G$ with $$\mathrm{sdeg}(B') = \mathrm{sdeg}(B)+2$$ and $|h_{\pm} B'| \leq |h_{\pm} B|$; repeated application of this results in a subset $A \in {\cal A}(G,m)$ with $|h_{\pm} A| \leq |h_{\pm} B|$, from which our result follows.

Let $$B=\{b_1,b_2,b_3,\dots,b_m\}$$ be an $m$-subset of $G$, and suppose that $-b_1 \not \in B$, $-b_2 \not \in B$, but $-b_3 \in B$.  Note that we may have $b_3=-b_3$; furthermore, the sets $\{\pm b_1\}$, $\{\pm b_2\}$, and $\{\pm b_3\}$ are pairwise disjoint.  Replacing $b_1$ by $-b_2$ in $B$, we let 
$$B'=\{-b_2, b_2, b_3,  \dots, b_m\}.$$  Then $B'$ has size $m$, and its degree of symmetry is exactly two more than that of $B$; we need to show that $|h_{\pm} B'| \leq |h_{\pm} B|$.  We shall, in fact, show that $h_{\pm} B' \subseteq h_{\pm} B$.  

By definition, $h_{\pm} B'$ is the collection of all elements of the form
$$g=\lambda_1(-b_2)+\lambda_2b_2+\lambda_3b_3+\cdots+\lambda_mb_m $$ where $\sum_{i=1}^m |\lambda_i|=h$.  Clearly, if $\lambda_1$ and $\lambda_2$ are of opposite sign or either one is zero, then $$|\lambda_2-\lambda_1|=|\lambda_1|+|\lambda_2|,$$ so 
$$g=(\lambda_2-\lambda_1)b_2+\lambda_3b_3+\cdots+\lambda_mb_m \in h_{\pm} B.$$ 

Suppose now that $\lambda_1$ and $\lambda_2$ are both positive; the case when they are both negative can be handled similarly.  Furthermore, we assume that $\lambda_1 \geq \lambda_2$; again, the reverse case is analogous.  

Assume first that $2b_3=0$; in this case we have $\lambda_3b_3=-\lambda_3b_3$, and thus we may assume that $\lambda_3 \geq 0$.  Observe that
$$g=(\lambda_2-\lambda_1)b_2+(2\lambda_1+\lambda_3)b_3+\lambda_4b_4+\cdots+\lambda_mb_m,$$ and
$$ |\lambda_2-\lambda_1|+|2\lambda_1+\lambda_3|+|\lambda_4|+\cdots+|\lambda_m|=h,$$
thus $g \in h_{\pm} B.$

Finally, suppose that $2b_3 \neq 0$; since $-b_3 \in B$, we must have $m \geq 4$, and without loss of generality we can assume that $b_4=-b_3$.  We can rewrite $g$ as follows:
$$g=\left\{
\begin{array}{ll}
(\lambda_2-\lambda_1)b_2+(\lambda_1+\lambda_3)b_3+(\lambda_1+\lambda_4)(-b_3)+\lambda_5b_5+\cdots+\lambda_mb_m & \mbox{if} \; \lambda_3 \geq 0, \lambda_4 \geq 0; \\ \\
(\lambda_2-\lambda_1)b_2+(\lambda_1+\lambda_3-\lambda_4)b_3+\lambda_1(-b_3)+\lambda_5b_5+\cdots+\lambda_mb_m & \mbox{if} \; \lambda_3 \geq 0, \lambda_4 \leq 0; \\ \\
(\lambda_2-\lambda_1)b_2+\lambda_1b_3+(\lambda_1-\lambda_3+\lambda_4)(-b_3)+\lambda_5b_5+\cdots+\lambda_mb_m & \mbox{if} \; \lambda_3 \leq 0, \lambda_4 \geq 0; \\ \\
(\lambda_2-\lambda_1)b_2+(\lambda_1-\lambda_4)b_3+(\lambda_1-\lambda_3)(-b_3)+\lambda_5b_5+\cdots+\lambda_mb_m & \mbox{if} \; \lambda_3 \leq 0, \lambda_4 \leq 0.
\end{array}
\right.$$
Since the expressions above show that $g \in h_{\pm} B$ in each case, our proof is complete.  $\Box$

\section{Cyclic groups}

In this section we prove that, when $G$ is cyclic,  then $\rho_{\pm} (G, m, h)$ agrees with $\rho (G, m, h)$ for all $m$ and $h$.
 
\begin{thm}  \label{cyclic} For all positive integers $n$, $m$, and $h$, we have 
$$\rho_{\pm} (\mathbb{Z}_n, m, h)= \rho (\mathbb{Z}_n, m, h).$$
\end{thm}

{\em Proof:}  Since the reverse inequality is obvious, it suffices to prove that $$\rho_{\pm} (\mathbb{Z}_n, m, h) \leq \rho (\mathbb{Z}_n, m, h).$$  Recall that $$\rho (\mathbb{Z}_n, m, h)=\min\{f_d (m,h) \; : \; d \in D(n)\}.$$  Observe that, for any symmetric subset $R$ of $G$ (that is, for every subset $R$ for which $R=-R$), we have $h_{\pm} R=hR$.
Our strategy is to find, for each $d \in D(n)$, a symmetric subset $R=R_d(n,m)$ of $\mathbb{Z}_n$ so that $|R| \geq m$ and $|hR| \leq f_d$; this will then imply that
$$\rho_{\pm} (\mathbb{Z}_n, m, h) \leq \min\{f_d (m,h) \; : \; d \in D(n)\}=\rho (\mathbb{Z}_n, m, h).$$  

We introduce some notations.  We write $n=2^an_0$, $d=2^bd_0$, and $\left \lceil m/d \right \rceil =2^cm_0$, where $a$, $b$, and $c$ are nonnegative integers and $n_0$, $d_0$, and $m_0$ are odd positive integers.  Our explicit construction of $R$ depends on whether $b+c \leq a$ or not.

Suppose first that $b+c \leq a$.  In this case, let $H$ be the subgroup of $G$ that has order $2^{c}d$, and set
$$R=\bigcup_{i=-\left \lfloor m_0/2 \right \rfloor}^{\left \lfloor m_0/2 \right \rfloor} (i+H).$$  Clearly, $R$ is symmetric; 
to see that $R$ has size at least $m$, note that for the index of $H$ in $G$ we have 
$$|G:H|=n/(2^cd)  \geq \left \lceil m/d \right \rceil/2^c =m_0=2 \left \lfloor m_0/2 \right \rfloor+1,$$hence $$|R|=\left( 2 \left \lfloor m_0/2 \right \rfloor+1 \right) \cdot |H|=d \left \lceil m/d \right \rceil \geq m.$$

To verify that $|hR| \leq f_d$, note that 
$$hR=\bigcup_{i=-h\left \lfloor m_0/2 \right \rfloor}^{h\left \lfloor m_0/2 \right \rfloor} (i+H),$$ so
\begin{eqnarray*}
|hR|&=& \min\{n,\left( 2h\left \lfloor m_0/2 \right \rfloor+1 \right) \cdot |H|\} \\
& \leq & \left( 2h\left \lfloor m_0/2 \right \rfloor+1 \right) \cdot |H|  \\
&=&(hm_0-h+1) \cdot 2^c d \\
&\leq &(2^c h m_0 -h+1)d \\
& = & f_d. 
\end{eqnarray*}

In the case when $b+c \geq a+1$, we let $H$ be the subgroup of $G$ that has order $2^ad_0$, and set
$$R=\bigcup_{i=1}^{2^{b+c-a-1}m_0} \left(\left \lfloor e/2 \right \rfloor +i +H \right) \cup \left(-\left \lfloor e/2 \right \rfloor -i +H\right),$$  where $e=n_0/d_0$. We see that $R$ is symmetric; in order to estimate $|R|$ and $|hR|$, we rewrite $R$ as follows.

Note that $e$ is an odd integer, and thus
$$  -\left \lfloor e/2 \right \rfloor = \left \lfloor e/2 \right \rfloor +1 - e;$$ furthermore, $e=n/|H|$ and thus $e$ is an element of $H$, and so 
$$-\left \lfloor e/2 \right \rfloor -i+H= \left \lfloor e/2 \right \rfloor +1-i+H$$ for every integer $i$.  
With this, we have
$$R=\bigcup_{i=-2^{b+c-a-1}m_0+1}^{2^{b+c-a-1}m_0} \left(\left \lfloor e/2 \right \rfloor +i +H \right).$$
To show that $R$ has size at least $m$, we see that, for the index of $H$ in $G$, we have  
$$|G:H|=n/(2^ad_0)  = 2^{b-a}n/d \geq 2^{b-a} \left \lceil m/d \right \rceil = 2^{b+c-a}m_0,$$hence $$|R|=\left( 2^{b+c-a}m_0 \right) \cdot |H|=d \left \lceil m/d \right \rceil \geq m.$$

Finally, $$hR=\bigcup_{i=-2^{b+c-a-1}hm_0+h}^{2^{b+c-a-1}hm_0} \left(H+ h \left \lfloor e/2 \right \rfloor +i \right),$$
so for $|hR|$ we get
\begin{eqnarray*}
|hR|&=& \min\{n,\left( 2^{b+c-a}hm_0-h+1 \right) \cdot |H|\} \\
& \leq & \left( 2^{b+c-a}hm_0-h+1 \right) \cdot |H|  \\
&=&\left( 2^{b+c-a}hm_0-h+1 \right) \cdot 2^ad_0 \\
&\leq &(2^c h m_0 -h+1)d \\
& = & f_d, 
\end{eqnarray*}
with which our proof is complete.  $\Box$

\section{Noncyclic groups}

Let us now turn to noncyclic groups.  We say that a finite abelian group $G$ has type $(n_1,\dots,n_r)$ if it is isomorphic to the invariant product  
$$\mathbb{Z}_{n_1} \times \cdots \times \mathbb{Z}_{n_r},$$ where $n_1 \geq 2$ and $n_i$ divides $n_{i+1}$ for each $i \in \{1,\dots, r-1\}$.  Here $r$ is the rank of $G$, $n_r$ is the exponent of $G$, and we still use the notation $n=\Pi_{i=1}^r n_i$ for the order of $G$. 

Recall that for the minimum size of the $h$-fold sumset of an $m$-subset of a group of order $n$ we have
$$\rho(G,m,h)=\min\{f_d (m,h) \; : \; d \in D(n)\}.$$  This, of course, implies that for signed sumsets we have the lower bound
$$\rho_{\pm} (G, m,h) \geq  \min \{f_d (m,h) \; : \; d \in D(n) \}.$$
It turns out that we can get an upper bound for $\rho_{\pm} (G, m,h)$ by minimizing $f_d$ for a certain subset of $D(n)$; 
more precisely, we establish the following result:

\begin{thm} \label{u pm with f}  The minimum size of the $h$-fold signed sumset of an $m$-subset of a group $G$ of type $(n_1,\dots,n_r)$ satisfies
$$\rho_{\pm} (G, m,h) \leq  \min \{f_d (m,h) \; : \; d \in D(G,m) \},$$
where $$D(G,m)=\{d \in D(n) \; : \; d= d_1 \cdots d_r, d_1 \in D(n_1), \dots, d_r \in D(n_r), dn_r \geq d_rm \}.$$
\end{thm}

Observe that, for cyclic groups of order $n$, $D(G,m)$ is simply $D(n)$.  

Theorem \ref{u pm with f} will be the immediate consequence of Propositions \ref{u pm is upper} and \ref{u pm with f prop} below.

\begin{prop}  \label{u pm is upper}
For every  group $G$ of type $(n_1,\dots,n_r)$ and order $n$, $m \leq n$, and $h \in \mathbb{N}$ we have
$$\rho_{\pm} \left(G, m, h   \right) \leq  u_{\pm} (G, m,h),$$ where
$$u_{\pm} (G, m,h)= \min \left \{ \Pi_{i=1}^r u(n_i,m_i,h) \; : \; m_1 \leq n_1, \dots, m_r \leq n_r, \Pi_{i=1}^r m_i \geq m \right\}.$$

\end{prop}

{\em Proof:}  For each $i=1,2,\dots,r$, let $m_i$ be an integer so that $m_i \leq n_i$ but $m_1 \cdots m_r \geq m$.  By Theorem \ref{cyclic}, for each $i$ we can find symmetric sets $A_i \subseteq \mathbb{Z}_{n_i}$ of size at least $m_i$ for which
$$|h_{\pm}A_i|=|hA_i|=u(n_i,m_i,h).$$  Therefore, $A_1 \times \cdots \times A_r$ is a symmetric subset of $Z_{n_1} \times \cdots \times Z_{n_r}$ of size at least $m_1 \cdots m_r$, so we have
\begin{eqnarray*}
\rho_{\pm} \left(\mathbb{Z}_{n_1} \times \cdots \times \mathbb{Z}_{n_r}, m, h   \right) &\leq & \rho_{\pm} \left(\mathbb{Z}_{n_1} \times \cdots \times \mathbb{Z}_{n_r}, m_1 \cdots m_r, h   \right) \\
& \leq & |h_{\pm} (A_1 \times \cdots \times A_r)| \\
& = & |h  (A_1 \times \cdots \times A_r)| \\
& \leq & |h  A_1 \times \cdots \times h A_r| \\
& = & u(n_1,m_1,h)  \cdots  u(n_r,m_r,h),
\end{eqnarray*} as claimed. $\Box$

\begin{prop} \label{u pm with f prop}  With the notations as introduced above, we have
$$u_{\pm} (G, m,h) = \min \{f_d (m,h) \; : \; d \in D(G,m) \}.$$

\end{prop}

{\em Proof:}  First, we prove that 
$$u_{\pm} (G, m,h) \leq \min \{f_d (m,h) \; : \; d \in D(G,m) \}.$$
Suppose that $d_1, \dots, d_r$ are positive integers so that $d_1 \in D(n_1), \dots, d_r \in D(n_r),$ and $d n_r \geq d_r m$.  Let $m_1=d_1, \dots, m_{r-1}=d_{r-1}$, and $m_r=\lceil d_r m/d \rceil$.  By assumption, $m_i \leq n_i$ for all $1 \leq i \leq r$, and we also have $m_1 \cdots m_r \geq m$; we will establish our claim by showing that $$u_{\pm} (G, m,h) \leq f_{d}(m, h).$$

Observe that, for each $1 \leq i \leq r-1$,
$$f_{d_i}(m_i,h)=f_{d_i}(d_i,h)=\left( h \left \lceil d_i/d_i \right \rceil -h+1 \right) d_i=d_i,$$
and
$$f_{d_r} (m_r,h)=f_{d_r} (\left \lceil d_r m/d \right \rceil,h)=\left( h \left \lceil \left \lceil d_rm/d \right \rceil /d_r \right \rceil -h+1 \right) d_r,$$ which, according to an identity for the ceiling function, equals
$$\left( h \left \lceil m/d \right \rceil -h+1 \right) d_r.$$
Therefore, $$f_{d_1}(m_1,h) \cdots f_{d_r} (m_r,h)  = \left( h \left \lceil m/d \right \rceil -h+1 \right) d=f_{d}(m, h).$$  Our claim now follows, since
$$u_{\pm} (G, m,h) \leq u(n_1,m_1,h) \cdots u(n_r,m_r,h) \leq f_{d_1}(m_1,h) \cdots f_{d_r} (m_r,h).$$

Conversely, we need to prove that 
\begin{eqnarray} \label{upper for r}
u_{\pm} (G, m,h) \geq \min \{f_d (m,h) \; : \; d \in D(G,m) \}.
\end{eqnarray}  As we have already mentioned, this holds for cyclic groups.  We will now prove that the inequality also holds for $r=2$; that is, for a group of type $(n_1,n_2)$ we have
\begin{eqnarray} \label{upper for r=2}
u_{\pm} (G, m,h) \geq \min \{f_{d_1d_2}(m,h) \; : \; d_1 \in D(n_1), d_2 \in D(n_2), d_1 n_2 \geq m \}.
\end{eqnarray}

Suppose that positive integers $m_1$ and $m_2$ are selected so that $m_1 \leq n_1$, $m_2 \leq n_2$, $m_1 m_2 \geq m$, and  
$$u_{\pm} (G, m,h)=u(n_1,m_1,h) \cdot u(n_2,m_2,h);$$ furthermore, suppose that integers $\delta_1$ and $\delta_2$ are chosen so that $\delta_1 \in D(n_1)$, $\delta_2 \in D(n_2)$, $u(n_1,m_1,h)=f_{\delta_1}(m_1,h)$, and $u(m_2,h)=f_{\delta_2}(m_2,h)$.  We need to prove that there are integers $d_1$ and $d_2$, so that $d_1 \in D(n_1)$, $d_2 \in D(n_2)$, $d_1 n_2 \geq m$, and 
\begin{eqnarray} \label{f_d_1d_2 leq f_delta1f_delta2} f_{d_1d_2}(m,h) &\leq& f_{\delta_1}(m_1,h) \cdot f_{\delta_2}(m_2,h).\end{eqnarray}
We will separate two cases depending on whether $\delta_1 n_2 \geq m$ or not.

In the case when $\delta_1 n_2 \geq m$, we show that $d_1=\delta_1$ and $d_2=\delta_2$ are appropriate choices.  Clearly, $d_1 \in D(n_1)$, $d_2 \in D(n_2)$, and $d_1 n_2 \geq m$, so we just need to show that $$f_{d_1d_2}(m,h) \leq f_{d_1}(m_1,h) \cdot f_{d_{2}}(m_2,h).$$ Since $m \leq m_1m_2$ and the function $f$ is nondecreasing in $m$, it suffices to prove that  
$$f_{d_1d_2}(m_1m_2,h) \leq f_{d_1}(m_1,h) \cdot f_{d_{2}}(m_2,h),$$ or, equivalently, that    
$$h \left \lceil (m_1m_2)/(d_1 d_2) \right \rceil -h+1  \leq \left( h \left \lceil m_1/d_1 \right \rceil -h+1 \right) \cdot \left( h \left \lceil m_2/d_2 \right \rceil -h+1 \right).$$  
Note that $$\lceil (m_1m_2)/(d_1 d_2) \rceil \leq \lceil m_1/d_1 \rceil \cdot \lceil m_2/d_2 \rceil,$$ so our inequality will follow once we prove that
$$h \left  \lceil m_1/d_1 \rceil \cdot \lceil m_2/d_2  \right \rceil -h+1  \leq \left( h \left \lceil m_1/d_1 \right \rceil -h+1 \right) \cdot \left( h \left \lceil m_2/d_2 \right \rceil -h+1 \right).$$ But this indeed holds as subtracting the left-hand side from the right-hand side yields
$$h(h-1)  \left( \left \lceil m_1/d_1 \right \rceil -1 \right) \left( \left \lceil m_2/d_2 \right \rceil -1 \right),$$ which is clearly nonnegative.

Suppose now that $\delta_1 n_2 < m$; we consider two subcases: when $m_2 \leq \delta_2$ and when $m_2 > \delta_2$.  

When  $\delta_1 n_2 < m$ and $m_2 \leq \delta_2$, we set $d_1=\gcd(n_1,\delta_2)$ and $d_2=\delta_1\delta_2/\gcd(n_1,\delta_2)$.  Then, clearly, $d_1 \in D(n_1)$; to see that $d_2 \in D(n_2)$, note that $n_1/d_1$ and $\delta_2/d_1$ are relatively prime integers that both divide $n_2/d_1$, so their product $n_1\delta_2/d_1^2$ divides $n_2/d_1$ as well, and therefore $n_1\delta_2/d_1$, and thus its divisor $d_2$, divide $n_2$.  Furthermore, since $n_1\delta_2/d_1$ divides $n_2$, we have $$d_1n_2 \geq n_1\delta_2 \geq m_1 m_2 \geq m.$$  It remains to be shown that (\ref{f_d_1d_2 leq f_delta1f_delta2}) holds, but since $d_1d_2=\delta_1\delta_2$, this follows as in the previous case.

Finally, suppose that $\delta_1 n_2 < m$ and $m_2 > \delta_2$; we now set $d_1=n_1$ and $d_2=\delta_1n_2/n_1$.   We see that $d_1 \in D(n_1)$,  $d_2 \in D(n_2)$, and $d_1n_2 \geq m$; we need to show that (\ref{f_d_1d_2 leq f_delta1f_delta2}) holds.

Let us denote $\left \lceil m_1/\delta_1 \right \rceil$ and $\left \lceil m_2/\delta_2 \right \rceil$  by $k_1$ and $k_2$, respectively; note that $m_2 > \delta_2$ implies that $k_2 \geq 2$, and $\delta_1 n_2 < m$ implies that $k_1 \geq 2$ as well, since
$$m_1 \geq m/m_2 > \delta_1 n_2/m_2 \geq \delta_1.$$  Therefore,
$$2(k_1-1)(k_2-1)=(k_1-2) (k_2-2)+(k_1k_2-2) \geq k_1k_2-2,$$ so, since $h \geq 2$, we get
$$h(h-1)(k_1-1) (k_2-1) \geq k_1k_2-2,$$ or, equivalently,
$$(hk_1-h+1) \cdot (hk_2-h+1) \geq (h+1)(k_1k_2-1).$$  Multiplying by $\delta_1\delta_2$ yields exactly
$$f_{\delta_1}(m_1,h) \cdot f_{\delta_2}(m_2,h)$$ on the left hand side; therefore, to prove (\ref{f_d_1d_2 leq f_delta1f_delta2}), we need to verify that
\begin{eqnarray} \label{f_d1d2 leq (h+1)}   f_{d_1d_2}(m,h)  \leq (h+1)(k_1k_2-1)\delta_1\delta_2.\end{eqnarray}
By definition,
$$f_{d_1d_2}(m,h)=f_{\delta_1n_2}(m,h)=\left( h\left \lceil m/(\delta_1n_2) \right \rceil -h+1    \right) \delta_1n_2.$$
But 
$$ \left \lceil \frac{m}{\delta_1n_2} \right \rceil \leq \left \lceil \frac{m_1m_2}{\delta_1n_2} \right \rceil \leq \left \lceil \frac{k_1k_2\delta_1\delta_2}{\delta_1n_2} \right \rceil = \left \lceil \frac{k_1k_2}{n_2/\delta_2} \right \rceil \leq \frac{k_1k_2+n_2/\delta_2-1}{n_2/\delta_2},$$
hence
\begin{eqnarray} \label{f_d1d2 leq again}  f_{d_1d_2}(m,h) \leq \left( h(k_1k_2-1)+n_2/\delta_2   \right) \delta_1\delta_2. \end{eqnarray}
Since we are under the assumption that $\delta_1 n_2 < m$, we have
$$\frac{n_2}{\delta_2} < \frac{m}{\delta_1\delta_2} \leq \frac{m_1m_2}{\delta_1\delta_2} \leq k_1k_2,$$ so the integer $n_2/\delta_2 $ can be at most $k_1k_2-1$, and thus (\ref{f_d1d2 leq again}) implies (\ref{f_d1d2 leq (h+1)}), completing the proof of (\ref{upper for r=2}).

In order to prove that (\ref{upper for r}) holds for any fixed $r>2$, we suppose that positive integers $m_1, \dots, m_r$ are selected so that $m_i \leq n_i$ for each $1 \leq i \leq r$, $m_1 \cdots m_r \geq m$, and  
$$u_{\pm} (G, m,h)=u(n_1,m_1,h) \cdots u(n_r,m_r,h).$$ Furthermore, we suppose that integers $\delta_1, \dots, \delta_r$ are chosen so that for each $1 \leq i \leq r$, $\delta_i \in D(n_i)$  and $u(n_i,m_i,h)=f_{\delta_i}(m_i,h)$.  
We will  prove that there are integers $d_1, \dots, d_r$, so that, for each $1 \leq i \leq r$, $d_i \in D(n_i)$, 
\begin{eqnarray} \label{d 1 d r-1}  d_1 \cdots d_{r-1} n_r &\geq & m, \end{eqnarray}
and 
\begin{eqnarray} \label{f_d_1d_r leq f_delta1f_deltar} f_{d_1 \cdots d_r}(m,h) &\leq& u_{\pm} (G, m,h) =f_{\delta_1}(m_1,h) \cdots f_{\delta_r}(m_r,h).\end{eqnarray}

We proceed by induction, and assume that (\ref{upper for r}) holds for $r-1$ terms and for $m'=m_2 \cdots m_r$; in particular, for a group G of rank $r-1$ and of type $(n_2,\dots,n_r)$ we have
$$u_{\pm} (G, m',h) \geq \min \{f_d (m',h) \; : \; d \in D(G, m') \}.$$  Therefore, we are able to find integers $\mu_2, \dots, \mu_r$ so that $\mu_i \in D(n_i)$ for each $2 \leq i \leq r$, 
\begin{eqnarray} \label{mu 1 mu r-1}  \mu_2 \cdots \mu_{r-1} n_r &\geq & m', \end{eqnarray}
and 
\begin{eqnarray}   \label{f_d_1d_r leq f_mu1f_mur} f_{\mu_2 \cdots \mu_r}(m',h) & \leq & u_{\pm} (G, m',h) \leq f_{\delta_2}(m_2,h) \cdots f_{\delta_r}(m_r,h).  \end{eqnarray} 

Furthermore, observing that by (\ref{mu 1 mu r-1}), $m''=\lceil m'/(\mu_2 \cdots \mu_{r-1}) \rceil $ is at most $n_r$,  from (\ref{upper for r=2}), for a group of rank 2 and of type $(n_1,n_r)$ we have 
$$u_{\pm} (G, m_1m'',h) \geq \min \{f_d (m_1m'',h) \; : \; d \in D(G,m_1m'') \},$$ and so there are integers $\nu_1 \in D(n_1)$ and $\nu_r \in D(n_r)$ for which 
\begin{eqnarray} \label{nu 1 nu r-1}  \nu_1 n_r &\geq & m_1m'', \end{eqnarray}
and 
\begin{eqnarray}   \label{f_d_1d_r leq f_nu1f_nur} f_{\nu_1\nu_r}(m_1m'',h) & \leq & u_{\pm} (G, m_1m'',h) \leq f_{\delta_1}(m_1,h) \cdot f_{\mu_r}(m'',h).  \end{eqnarray}

Now let $d_1=\nu_1$, $d_r=\nu_r$, and $d_i=\mu_i$ for $2 \leq i \leq r-1$.  We immediately see that, with these notations, (\ref{d 1 d r-1})  holds, since, by (\ref{nu 1 nu r-1}),
$$d_1 \cdots d_{r-1}n_r = \nu_1 \mu_2 \cdots \mu_{r-1} n_r \geq m_1 \mu_2 \cdots \mu_{r-1} m'' \geq m_1 m' = m_1 \cdots m_r \geq m.$$

To see that (\ref{f_d_1d_r leq f_delta1f_deltar}) holds, note that, for the left-hand side we have
\begin{eqnarray*}
f_{d_1 \cdots d_r}(m,h) &=& f_{\nu_1\nu_r\mu_2 \cdots \mu_{r-1}}(m,h) \\
   & \leq & f_{\nu_1\nu_r\mu_2 \cdots \mu_{r-1}}(m_1 m'' \mu_2 \cdots \mu_{r-1},h) \\
& = & \left(  h \left \lceil (m_1m'')/(\nu_1 \nu_r) \right \rceil -h+1 \right) \nu_1\nu_r\mu_2 \cdots \mu_{r-1} \\
  & = & f_{\nu_1 \nu_r} (m_1m'',h) \mu_2 \cdots \mu_{r-1};
\end{eqnarray*}
and, for the right-hand side of (\ref{f_d_1d_r leq f_delta1f_deltar}), using (\ref{f_d_1d_r leq f_mu1f_mur}), we see that
\begin{eqnarray*}
f_{\delta_1}(m_1,h) \cdots f_{\delta_r}(m_r,h) & \geq & f_{\delta_1}(m_1,h)  f_{\mu_2 \cdots \mu_r}(m',h) \\
&=& f_{\delta_1}(m_1,h)  \left(  h \left \lceil m'/(\mu_2 \cdots \mu_r) \right \rceil -h+1 \right)\mu_2 \cdots \mu_r \\
&=& f_{\delta_1}(m_1,h)  \left(  h \left \lceil m''/\mu_r \right \rceil -h+1 \right)\mu_2 \cdots \mu_r \\
&=& f_{\delta_1}(m_1,h)  f_{\mu_r}(m'',h)\mu_2 \cdots \mu_r.
\end{eqnarray*}
Therefore, (\ref{f_d_1d_r leq f_delta1f_deltar}) follows from (\ref{f_d_1d_r leq f_nu1f_nur}).  With this, the proof of (\ref{upper for r}), and thus of Proposition \ref{u pm with f prop}, is complete.
$\Box$

Our next result exhibits a situation where the upper bound of Proposition \ref{u pm is upper}, and thus of Theorem \ref{u pm with f}, is not tight:

\begin{prop}
If $G$ is a noncyclic group of odd order $n$ and type $(n_1, \dots, n_r)$, then
$$\rho_{\pm} \left(G, (n-1)/2, 2   \right) \leq n-1,$$ but $$u_{\pm} (G, (n-1)/2,2)=n.$$ 

\end{prop}

{\em Proof:}  Note that every element of $G \setminus \{0\}$ has order at least 3, thus there is a subset $A$ of $G \setminus \{0\}$ with which $G \setminus \{0\}$ can be partitioned into $A$ and $-A$.  Since $|A|=(n-1)/2$ and  $0 \not \in 2_{\pm}A$, we have  $$\rho_{\pm} \left(G, (n-1)/2, 2   \right) \leq n-1.$$ 

To prove our second claim, note that for each $i \in \{1,\dots,r\}$, 
$$n/n_i \cdot (n_i-1)/2 < (n-1)/2.$$  Therefore, if positive integers $m_1, \dots, m_r$ satisfy $m_i \leq n_i$ for each $i \in \{1,\dots,r\}$ and $$m_1\cdots m_r \geq (n-1)/2,$$ then  we must have $m_i \geq (n_i+1)/2$, and thus $u(n_i,m_i,2)=n_i$, for each $i \in \{1,\dots,r\}$, from which our claim follows.  $\Box$

A bit more generally, if $d$ is an odd element of $D(n)$ so that $d \geq 2m+1$, then the same argument yields 
$$\rho_{\pm} \left(G, m, 2   \right) \leq d-1,$$
and therefore we have the following:

\begin{cor}
Suppose that $G$ is an abelian group of order $n$ and type $(n_1, \dots, n_r)$.  Let $m \leq n$, and let $d_m$ be the smallest odd element of $D(n)$ that is at least $2m+1$; if no such element exists, set $d_m=\infty$.  We then have
$$\rho_{\pm} \left(G, m, 2   \right) \leq \min\{u_{\pm}(G,m,2), d_m-1\}.$$

\end{cor}

We are not aware of any subsets with smaller signed sumset size, and we believe that the following holds:

\begin{conj} \label{conj for rho pm}
Suppose that $G$ is an abelian  group of order $n$ and type $(n_1, \dots, n_r)$.  Let $m \leq n$ and $h \geq 2$.  

If $h \geq 3$, then $$\rho_{\pm} \left(G, m, h   \right) =u_{\pm}(G,m,h).$$

If each odd divisor of $n$ is less than $2m$, then $$\rho_{\pm} \left(G, m, 2  \right) =u_{\pm}(G,m,2).$$

If there are odd divisors of $n$ greater than $2m$, let $d_m$ be the smallest one.  We then have
$$\rho_{\pm} \left(G, m, 2   \right) = \min\{u_{\pm}(G,m,2), d_m-1\}.$$
\end{conj}

\section{An example}

Trivially, if $G$ is an elementary abelian 2-group, then $\rho_{\pm} \left(G, m, h   \right)$ agrees with $\rho \left(G, m, h   \right)$, and it is not hard to see that this is also true if $G$ is any 2-group.  More generally still, as an application to Theorem \ref{u pm with f}, we prove the following:

\begin{prop} \label{example}
If there is no odd prime $p$ for which $\mathbb{Z}_p^2$ is isomorphic to a subgroup of $G$, then  
$$\rho_{\pm} \left(G, m, h   \right) = \rho \left(G, m, h   \right).$$

\end{prop} 

{\em Proof:}  Suppose that $G$ is of order $n$ and of type $(n_1,\dots,n_r)$; by Theorem \ref{cyclic}, we may assume that $r \geq 2$.  

Let $d \in D(n)$ be such that $$\rho \left(G, m, h   \right)=u(n,m,h)=f_d(m,h).$$  By Theorem \ref{u pm with f}, it suffices to prove that $d \in D(G,m)$.  

Our assumption that there is no odd prime $p$ for which $\mathbb{Z}_p^2$ is isomorphic to a subgroup of $G$ is equivalent to saying that $n_1 \cdots n_{r-1}$ is a power of 2; let $$n_1 \cdots n_{r-1}=2^{k_1}.$$ Furthermore, we write $$n_r=2^{k_2} \cdot c_2$$ and $$d=2^{k_3} \cdot c_3,$$ where $k_2$ and $k_3$ are nonnegative integers, and $c_2$ and $c_3$ are odd.  Note that 
\begin{eqnarray} \label{k1 and k2 and k3}
k_1+k_2 & \geq & k_3,
\end{eqnarray} and $c_2$ must be divisible by $c_3$. 

Now if $m \leq n_r$, then clearly $d \in D(G,m)$, so assume that $m \geq n_r+1$, and thus there is a nonnegative integer $k$ for which
$$2^k \cdot n_r +1 \leq m \leq 2^{k+1} \cdot n_r.$$  Note that we must then have 
\begin{eqnarray} \label{k1 and k}
k_1 & \geq & k+1.
\end{eqnarray}

We claim that we also have 
\begin{eqnarray} \label{k3 and k2 and k}
k_3 &\geq & k_2+k+1.
\end{eqnarray}  Indeed, 
\begin{eqnarray*}
u(n,m,h) & = & f_d(m,h) \\ 
 & = & \left(h \cdot  \left \lceil m/h \right \rceil-h +1 \right) \cdot d \\ 
& \geq & \left( h \cdot \left \lceil \frac{2^k \cdot n_r + 1}{d} \right \rceil-h +1 \right) \cdot d.
\end{eqnarray*}
On the other hand, from (\ref{k1 and k}) we see that $G$ contains a subgroup of order $2^{k+1} \cdot n_r$, and thus
\begin{eqnarray*}
u(n,m,h) & \leq & 2^{k+1} \cdot n_r \\ 
 & < & h \cdot 2^{k} \cdot n_r +d \\ 
& = & \left( h \cdot \frac{2^k \cdot n_r + d}{d} -h +1 \right) \cdot d.
\end{eqnarray*}
Therefore, 
$$ \left \lceil \frac{2^k \cdot n_r + 1}{d} \right \rceil < \frac{2^k \cdot n_r + d}{d},$$
which yields that $2^k \cdot n_r$ cannot be divisible by $d$, that is, $2^{k+k_2} \cdot c_2$ cannot be divisible by $2^{k_3} \cdot c_3$, proving (\ref{k3 and k2 and k}). 

Now let $$d_r=2^{k_2} \cdot c_3.$$  Then $d_r$ is a divisor of $n_r$; furthermore, by (\ref{k3 and k2 and k}), $d/d_r=2^{k_3-k_2}$ is an integer, and by (\ref{k1 and k2 and k3}), it is a divisor of $n_1 \cdots n_{r-1}$.  Using (\ref{k3 and k2 and k}) again, we have
$$d \cdot n_r = 2^{k_3} \cdot c_3 \cdot n_r \geq 2^{k_2+k+1} \cdot c_3 \cdot n_r = d_r \cdot 2^{k+1} \cdot n_r \geq d_r \cdot m,$$ so $d \in D(G,m)$, as claimed.  $\Box$

Having a subgroup that is isomorphic to $\mathbb{Z}_p^2$ for an odd prime $p$ is thus a necessary condition for $\rho_{\pm} \left(G, m, h   \right)$ to be greater than $\rho \left(G, m, h   \right)$.  We study $\mathbb{Z}_p^2$, and, more generally, elementary abelian groups, in the upcoming paper \cite{BajMat:2014b}.


\begin{thebibliography}{99}


\bibitem{Baj:2000a} B. Bajnok, Spherical Designs and Generalized Sum-Free Sets in Abelian Groups.  Special issue dedicated to Dr. Jaap Seidel on the occasion of his 80th birthday (Oisterwijk, 1999).  \emph{Des. Codes Cryptogr.}  {\bf 21} (2000), no. 1--3, 11-18.

\bibitem{Baj:2004a} B. Bajnok, The Spanning Number and the Independence Number of a Subset of an Abelian Group.  In \emph{Number Theory,} D. Chudnovsky, G. Chudnovsky, and M. Nathalson (Ed.), Springer-Verlag (2004), 1-16.

\bibitem{BajMat:2014b} B. Bajnok and R. Matzke, On the Minimum Size of Signed Sumsets in Elementary Abelian Groups, www.arxiv.org (2014).

\bibitem{BajRuz:2003a} B. Bajnok and I. Ruzsa, The Independence Number of a Subset of an Abelian Group.  \emph{Integers} {\bf 3} (2003), Paper No. A2, 23 pp.

\bibitem{Cau:1813a} A.-L. Cauchy, Recherches sur les nombres, {\em J. \'Ecole Polytechnique} {\bf 9} (1813), 99--123.

\bibitem{Dav:1935a} H. Davenport, On the addition of residue classes, {\em J. London Math. Soc.} {\bf 10} (1935), 30--32.

\bibitem{Dav:1947a} H. Davenport, A historical note, {\em J. London Math. Soc.} {\bf 22} (1947), 100--101.


\bibitem{EliKer:2005a} S. Eliahou and M. Kervaire, Old and new formulas for the Hopf--Stiefel and related functions, {\em Expo. Math.}, {\bf 23} (2005), no. 2, 127--145.

\bibitem{EliKer:2007a} S. Eliahou and M. Kervaire, Some extensions of the Cauchy--Davenport Theorem, {\em Electron. Notes in Discrete Math.}, {\bf 28} (2007) 557--564.


\bibitem{EliKerPla:2003a} S. Eliahou, M. Kervaire, and A. Plagne, Optimally small sumsets in finite abelian groups, {\em J. Number Theory}, {\bf 101} (2003), 338--348.


\bibitem{Kar:2006a} Gy. K\'arolyi, A note on the Hopf--Stiefel function.  {\em European J. Combin.}, {\bf 27} (2006), 1135--1137.


\bibitem{KloLev:2003a} B. Klopsch and V. F. Lev, How long does it take to generate a group?  {\em J. Algebra}, {\bf 261} (2003), 145--171.

\bibitem{KloLev:2009a} B. Klopsch and V. F. Lev, Generating abelian groups by addition only.  {\em Forum Math.}, {\bf 21} (2009), no. 1, 23--41.


\bibitem{Pla:2003a} A. Plagne, Additive number theory sheds extra light on the Hopf--Stiefel $\circ$ function, {\em Enseign. Math., II S\'er}, {\bf 49}(2003), no. 1--2, 109--116.

\bibitem{Pla:2006a} A. Plagne, Optimally small sumsets in groups, I. The supersmall sumset property, the $\mu_G^{(k)}$ and the $\nu_G^{(k)}$ functions, {\em Unif. Distrib. Theory}, {\bf 1} (2006), no. 1, 27--44.

\bibitem{Sha:1984a} D. Shapiro, Products of sums of squares, {\em Expo. Math.}, {\bf 2} (1984), 235--261.




\end{thebibliography}
\end{document}